\documentclass[12pt]{amsart}
\usepackage[colorlinks=true,citecolor=black,linkcolor=black,urlcolor=blue]{hyperref}
\usepackage{amsmath}
\usepackage[english,  activeacute]{babel}
\usepackage[utf8]{inputenc}
\usepackage{amssymb}
\usepackage{amsthm}
\usepackage{graphics,graphicx}
\usepackage{array}
\usepackage{bm}
\usepackage{a4wide}
\usepackage{color, url}
\usepackage{float}
\usepackage{pgf, tikz}
\usepackage{amsmath,tabularx}
\usepackage[T1]{fontenc}

\usepackage{algorithm}
\usepackage[noend]{algpseudocode}

\theoremstyle{plain}
\newtheorem{theorem}{Theorem}[section]
\newtheorem{proposition}[theorem]{Proposition}
\newtheorem{corollary}[theorem]{Corollary}

\theoremstyle{definition}
\newtheorem{definition}[theorem]{Definition}
\newtheorem{example}[theorem]{Example}
\newtheorem{remark}[theorem]{Remark}

\def\sttf2#1#2{\left[\!\!\left[#1\atop#2\right]\!\!\right]}

\RequirePackage{fix-cm}


\begin{document}

\title{Over-Mahonian numbers: Basic properties and unimodality}
\author{Ali Kessouri}
\address{\noindent Department of Mathematics, University of Ferhat Abbas
Setif 1, Algeria}
\email{ali.kessouri@ensta.edu.dz or ali.kessouri@univ-setif.dz}
\author{Moussa Ahmia}
\address{\noindent Department of Mathematics, University of Mohamed Seddik
Benyahia, LMAM laboratory, BP 98 Ouled Aissa Jijel 18000, Algeria}
\email{moussa.ahmia@univ-jijel.dz or ahmiamoussa@gmail.com}
\author{Salim Mesbahi}
\address{\noindent Department of Mathematics, University of Ferhat Abbas
Setif 1, Algeria}
\email{salim.mesbahi@univ-setif.dz}
\date{\today }

\begin{abstract}
In this paper, we introduce the concept of the over-Mahonian number, which counts the overlined permutations of length $n$ with $k$ inversions, allowing the first elements associated with the inversions to be independently overlined or not. We explore its properties and combinatorial interpretations through lattice paths, overpartitions, and tilings, and provide a combinatorial proof demonstrating that these numbers form a log-concave and unimodal sequence.
\end{abstract}

\maketitle

\section{Introduction}
A permutation $\pi=\left(
  \begin{array}{ccc}
    1 & \cdots & n \\
    \pi(1) & \cdots & \pi(n) \\
  \end{array}\right) $ is an arrangement or reordering of a set of elements. For additional information on the permutations, refer to \cite{bona, Sta86}.

\medskip

In mathematics, permutation statistics are functions that assign numerical values to permutations, often capturing combinatorial or algebraic properties.  Some of the most studied statistics include:
\begin{description}
\item[$\bullet$] {\bf Inversions:} The number of pairs $(i,j)$ where $i<j$ but $\pi(i)>\pi(j)$. This measures how "out of order" a permutation is.
\item[$\bullet$] {\bf Descent:} The number of positions $i$ such that $\pi(i)>\pi(i+1)$. A descent represents a local drop in the value sequence of the permutation.
\item[$\bullet$] {\bf Major Index:} The sum of all descent positions. For a permutation $\pi$, it is defined as $maj(\pi)=\sum_{i \text{\ is a descent}}i$.
\item[$\bullet$] {\bf Fixed Points:} The number of indices $i$ such that $\pi(i)=i$.
\item[$\bullet$] {\bf Cycle Structure:} Permutations can be decomposed into disjoint cycles. The number and lengths of these cycles are key permutation statistics, with implications in algebra and combinatorics.
\item[$\bullet$] {\bf Excedance:}  The number of indices $i$ such that $\pi(i)>i$.   
\end{description}
\noindent These statistics have significant applications in algebra, geometry, and computer science, especially in studying the symmetric group $S_n$, which consists of all permutations of a set with $n$ elements.

\medskip

The enumeration of permutations of length $n$ by their number of inversions, along with the study of $i(n,k)$, of permutations of length $n$ with $k$ inversions, is a foundational topic in combinatorics. The most famous result is given by the following equation \cite{Ro}:
\begin{equation}\label{giv}
\sum_{k=0}^{\binom{n}{2}}i(n,k)x^k=(1+x)\cdots(1+x+\cdots+x^{n-1}),
\end{equation}
where $i(n,k)$ represents the Mahonian number. MacMahon \cite{mac} established that this number corresponds to the count of permutations of length $n$ with a major index of $k$.

\medskip

Mahonian numbers can be extended to other finite reflection groups beyond the symmetric group. For example:
\begin{description}
  \item[$\bullet$] {\bf Coxeter groups:} Permutations in these groups can have statistics analogous to inversions or major indices.
  \item[$\bullet$]{\bf Signed permutations:} In type $B$ (also called hyperoctahedral group $B_n$) and type $D$ Coxeter groups, corresponding Mahonian statistics are defined.

\end{description}
\noindent For a detailed exploration of the combinatorics of the final groups, we direct readers to Björner and Brenti's book \cite{Bren}.
\medskip

The recent generalization of Mahonian numbers through generalized symmetric groups has sparked significant interest among researchers specializing in combinatorics. This generalization is rooted in both classical and contemporary permutation statistics. For instance, the Mahonian numbers of type $B$ extend the classical Mahonian numbers, which are associated with the symmetric group $S_n$, to the {\bf hyperoctahedral group $B_n$} using inversions of type $B$: inversions that include both pairs and signs. The group $B_n$ consists of signed permutations of $n$ elements, where each element can independently be positive or negative. For more details about the Mahonian numbers of type $B$, see \cite{Ar,KAAM}. 

\medskip

If we consider overlined permutations, where each element can independently be overlined or not, rather than signed permutations, we remain entirely within the hyperoctahedral group $B_n$. The hyperoctahedral groups have been extensively studied, as evidenced by works such as  \cite{r1, r2, r5, r6, r7, r11}.

\medskip

In our paper, we aim to introduce the concept of an overlined permutation $\sigma$ of $n$ elements, where the first elements associated with the inversions can be independently overlined or no. In such a permutation, each element $\sigma(i)$, for $1\leq i\leq n-1$, can independently be overlined or not, subject to the condition that there exists an index $j>i$ such that $\sigma(i)>\sigma(j)$. These permutations form a subgroup of the hyperoctahedral group $B_n$, which we denote by  $B'_n$. We focus on counting such permutations with exactly $k$ inversions, denoted by $i_{B'}(n, k)$, and refer to these counts as {\bf over-Mahonian numbers}. Furthermore, we study their combinatorial interpretations, identities, and provide combinatorial proofs of their log-concavity and unimodality.    

\bigskip

The paper is structured into four sections. In Section \ref{sec:2}, we define the over-Mahonian numbers as the counts of overlined permutations of length $n$ with $k$ inversions, where the first elements associated with the inversions can be independently overlined or not, within the hyperoctahedral group $B_n$, and provide a combinatorial proof of several recurrence relations and identities associated with these numbers. We also present their generating function and derive a key identity involving the double factorial. Section \ref{sec:3} offers combinatorial interpretations of the over-Mahonian numbers through lattice paths, overpartitions, and tilings. In Section \ref{sec:4}, we prove combinatorially that these numbers form a log-concave and thus unimodal sequence, using an appropriate injection. Finally, in the fourth section, we pose a question regarding the number and positions of the modes in the sequence of over-Mahonian numbers.
\section{Basic properties of over-Mahonian numbers}\label{sec:2}
Let $[n]$ denote the set $\{1,2,\ldots,n\}$, and let $\sigma=\sigma_1\sigma_2\cdots\sigma_n$ represents a permutation of $[n]$ of length $n$. The set of all such permutations forms the symmetric group $S_n$, whose identity element is the permutation  $\iota=\iota_1\iota_2\cdots\iota_n$, defined by $\iota_i=i$ for all $1\leq i\leq n$.  

\begin{definition} The backward  (or the inverse) permutation of $\sigma$  is the  permutation $\sigma'$ defined by $\sigma'_i=\sigma_{n+1-i}$ for  all $1\leq i\leq n$.
\end{definition}

\begin{remark}
The only permutation with no inversion is the identity permutation $\iota$, while the permutation of length $n$  with the maximum number of inversions is $\alpha=\alpha_1\alpha_2\cdots\alpha_n$, where $\alpha_i=n+1-i$, which has $\binom{n}{2}$ inversions;  furthermore, a permutation and its inverse have the same number of inversions.
\end{remark}

\begin{definition} A pair $(\sigma_i,\sigma_j)$  is termed a backward inversion of the permutation $\sigma$  if $i<j$ and $\sigma_i<\sigma_j$.
\end{definition} 

\begin{remark} For a permutation of length $n$ with $k$ inversions, the number of backward inversions is given by $\binom{n}{2}-k$ \cite{dex}.
\end{remark}

Consider the set of $2n$ symbols 
$$\Sigma_{2,n}:=\left\{1,\ldots,n,\overline{1},\ldots,\overline{n}\right\}.$$

An element denoted as $\overline{i}$ is referred to as an overlined element. An overlined permutation $\pi$ is a permutation defined on the set $\Sigma_{2,n}$ that satisfies the property $\pi(\overline{a}) = \overline{\pi(a)}$ for all $a \in \Sigma_{2,n}$. For instance, the following is an example of an overlined permutation of $\Sigma_{2,4}$:
$$\pi=\left(
  \begin{array}{cccc}
    1 & 2 & 3 & 4\\
    \overline{4} & \overline{2} & 3& 1\\
  \end{array}\right).$$
By omitting the first row, we obtain the one-line notation $\overline{4}\text{\ }\overline{2}3 1$. Let $G_{2,n}$ denote the set of all overlined permutations of $\Sigma_{2,n}$. This set, $G_{2,n}$, corresponds to the hyperoctahedral group $B_n$, also known as the Coxeter group of type $B$, and has cardinality $|G_{2,n}| = 2^n n!$. In algebraic combinatorics, $G_{2,n}$ is identified as the wreath product $C_2 \wr S_n$, combining the symmetric group on $[n]$ with the cyclic group $C_2$ on $\{0,1\}$. However, the group structure of $G_{2,n}$ is not directly relevant to this work.

\medskip
In this work, we aim to define a subset of the set of overlined permutations $G_{2,n}$, by imposing specific restrictions, then enumerate its elements while accounting for the number of inversions. This subset, denoted by $B'_n\subset G_{2,n}$ consists of overlined permutations where the first elements associated with inversions can independently be overlined or not. For example, for $n=3$, we have
\[B'_3=\{123,132, 1\overline{3}2,213,\overline{2}13,231,\overline{2}31,2\overline{3}1,\overline{2}\text{\ }\overline{3}1,312,\overline{3}12,321,\overline{3}21,3\overline{2}1,\overline{3}\text{\ }\overline{2}1\}.\]
A permutation $\sigma$ on $S_n$ is an involution if $\sigma^2(i)=i$ for all $i=1,\ldots,n$, and it is {\bf fixed-point-free} if no element is mapped to itself, meaning $\sigma(i)\neq  i$ for all $i$. Thus, the cardinality $|B'_n|$ is equal to the number of fixed-point-free involutions in the symmetric group $S_{2n}$. $|B'_n|$ is equal also to the number of permutations in the symmetric group $S_{2n}$ whose cycle decomposition is a product of $n$ disjoint transpositions.

Based on the definition of $B'_n$, in the next, we introduce the concept of over-Mahonian numbers along with their fundamental properties.
\begin{definition}\label{df} Let $i_{B'}(n,k)$ represent the number of overlined permutations in $B'_n$ of length $n$ that contain exactly $k$ inversions. We refer to this quantity as the \textbf{over-Mahonian number}. 
\end{definition}
\begin{example}
Consider the overlined permutations of $B'_3$ with 2 inversions: $$312, \mathbf{\overline{3}}12, 231, \mathbf{\overline{2}}31, 2\mathbf{\overline{3}}1, \mathbf{\overline{2}}\hspace{0.05cm}\mathbf{\overline{3}}1.$$
Therefore, $i_{B'}(3,2)=6$.
\end{example}

From Definition \ref{df}, we can derive the recurrence relation for the over-Mahonian number $i_{B'}(n,k)$ as follows:
\begin{theorem}\label{thmm}
For positive integers $n$ and $k$, the over-Mahonian number satisfies the recurrence relation
\begin{equation} \label{eq1} 
i_{B'}(n,k) = i_{B'}(n-1,k) + 2 \sum_{j=1}^{n-1} i_{B'}(n-1,k-j), 
\end{equation} 
with the initial conditions $i_{B'}(n,0) = 1$ and $i_{B'}(n,k) = 0$ unless $0 \leq k \leq \binom{n}{2}$. 
\end{theorem}
\begin{proof}
We prove this theorem combinatorially using the inversion combinatorial interpretation given in Definition \ref{df}.

Consider a permutation $\sigma=\sigma_1\cdots \sigma_n$ of length $n$ with $k$ overlined inversions. Removing the first element $\sigma_1$ form $\sigma$ results in either: 

\begin{enumerate}
  \item A permutation of length $n-1$ with $k$ overlined inversions, if  $\sigma_1$ does not form an overlined inversion with any other element. This corresponds to $i_{B'}(n-1,k$, or
      
  \item A permutation of length $n-1$ with $k-j$ overlined inversions, if $\sigma_1$ either forms  $j$ overlined inversions (where $1\leq j\leq n-1$) with other elements or is itself overlined and forms $j$ overlined inversions with the others. In both cases, the corresponding interpretation is $i_{B'}(n-1,k-j)$.
\end{enumerate}
Summing these possibilities yields $2 \sum_{j=1}^{n-1}i_{B'}(n-1,k-j)$, which leads to the equality \eqref{eq1}.
\end{proof}

\medskip

From Theorem \ref{thmm}, we derive the following recurrence relation for the over-Mahonian numbers $i_{B'}(n,k)$:
\begin{proposition}For $0\leq k \leq \binom{n}{2}$, the over-Mahonian numbers satisfy the recurrence 
\begin{equation}\label{rcc}
i_{B'}(n,k)=i_{B'}(n,k-1)+i_{B'}(n-1,k)+i_{B'}(n-1,k-1)-2i_{B'}(n-1,k-n).
\end{equation}
\end{proposition}
\begin{proof}
From relation \eqref{eq1}, we have the following expressions:
\begin{equation} \label{eq11}
i_{B'}(n,k)=i_{B'}(n-1,k)+2i_{B'}(n-1,k-1)+\cdots +2i_{B'}(n-1,k-n+1)
\end{equation}
and 
\begin{equation} \label{eq12}
i_{B'}(n,k-1)=i_{B'}(n-1,k-1)+2i_{B'}(n-1,k-2)+\cdots +2i_{B'}(n-1,k-n). 
\end{equation}
Using relation \eqref{eq12}, we substitute into the equation for $i_{B'}(n,k)$ : 
\begin{align*} 
i_{B'}(n,k) &= i_{B'}(n-1,k) + i_{B'}(n-1,k-1) +  [i_{B'}(n-1,k-1) + 2 i_{B'}(n-1,k-2) \\ 
&\quad +\cdots + 2 i_{B'}(n-1,k-n)]- 2  i_{B'}(n-1,k-n) \\
 &= i_{B'}(n-1,k) + i_{B'}(n-1,k-1) + i_{B'}(n,k-1) - 2 i_{B'}(n-1,k-n) \\
 &= i_{B'}(n-1,k) + i_{B'}(n-1,k) + i_{B'}(n-1,k-1) - 2 i_{B'}(n-1,k-n). 
\end{align*} This completes the proof.
\end{proof}

Based on Definition \ref{df}, we can also derive the following identities.
\begin{proposition}For any integer $n\geq 1$, the following hold:
\begin{enumerate}
  \item- $i_{B'}(n,\binom{n}{2})=2^{n-1}$,
 \item- $i_{B'}(n,1)=2\times(n-1)$,
\item- $i_{B'}(n,k)$ is even for $n\geq 2$ and $k\geq 1$.
\end{enumerate}
\end{proposition}
\begin{proof}
\begin{enumerate}
\item The expression $i_{B'}(n,\binom{n}{2})$ counts the number of permutations of length $n$ that have $\binom{n}{2}$ overlined inversions. This implies that the permutation is in the form $\sigma_1 > \sigma_2 > \cdots > \sigma_{n-1} > \sigma_n$, where the entries $\sigma_1, \ldots, \sigma_{n-1}$ may or may not be overlined. Thus, there are $2^{n-1}$ possible permutations.
\item The expression $i_{B'}(n,1)$ counts the number of permutations of length $n$ that have exactly one overlined inversion. This corresponds to a permutation in the form $\sigma_1 < \sigma_2 < \cdots < \sigma_j > \sigma_{j+1} < \cdots < \sigma_n$, where $\sigma_j$ may or may not be overlined, for each $1 \leq j \leq (n-1)$. Therefore, there are $2 \times (n-1)$ possible permutations.
\item For $n \geq 2$ and $k \geq 1$, we can have overlined inversions. For each inversion $(\sigma_i, \sigma_j)$, there are two possibilities: either $\sigma_i$ is overlined or it is not. Consequently, the total number of permutations of length $n$ with exactly $k$ overlined inversions is even.   
\end{enumerate}    
\end{proof}

The over-Mahonian numbers $i_{B'}(n,k)$ satisfy the following row-generating function:
\begin{theorem}\label{thm2}For any positive integer $n$,
\begin{equation*}
\sum_{k=0}^{\binom{n}{2}}i_{B'}(n,k)z^{k}=%
\begin{cases}
1, & \text{if \ \ } n=1, \\ 
(1+2z)\cdots(1+2z+\cdots+2z^{n-1}), & \text{if \ \ }n>1.%
\end{cases}%
\end{equation*}
\end{theorem}
\begin{proof}
Let $f_{n}(z)=\sum_{k=0}^{\binom{n}{2}}i_{B'}(n,k)z^{k}$ be the row generating function of the over-Mahonian number.

\medskip

For $n=1$, the result is trivial. For $n\geq 2$, using equation \eqref{eq1}, we have
 \begin{align*}
   f_{n}(z) & =\sum_{k=0}^{\binom{n}{2}}i_{B'}(n,k)z^{k}= \sum_{k=0}^{\binom{n}{2}}\left( i_{B'}(n-1,k)+2\sum_{j=1} ^{n-1}i_{B'}(n-1,k-j)\right) z^{k}.
  \end{align*}
  This expands to
 \begin{align*} 
   f_{n}(z) & = \sum_{k=0}^{\binom{n}{2}}i_{B'}(n-1,k)z^{k}+2\sum_{k=0}^{\binom{n}{2}}\sum_{j=1}^{n-1}i_{B'}(n-1,k-j) z^{k}.
   \end{align*}    
 The first term is simply $ f_{n-1}(z)$, and the second term can be written as
 \begin{align*}  
 2\sum_{j=1}^{n-1}z^{j}\sum_{k=j} ^{\binom{n}{2}-j}i_{B'}(n-1,k-j)z^{k-j}&=2\sum_{j=1}^{n-1}z^{j}\sum_{k=0}^{\binom{n-1}{2}}i_{B'}(n-1,k)z^{k}.
 \end{align*}
Thus, we have
 \begin{align*} 
    f_{n}(z)&= f_{n-1}(z)\left(1+2\sum_{j=1}^{n-1}z^j\right). 
 \end{align*} 

Iterating this process, we obtain the final result:
$$f_{n}(z)=(1+2z)\cdots(1+2z+\cdots+2z^{n-1}).$$ 
\end{proof}
By setting $z=1$ in Theorem \ref{thm2}, we obtain the following interesting result:
\begin{corollary}
For any positive integer $n$, the over-Mahonian numbers satisfy the identity: 
\begin{equation*} |B'_n|=\sum_{k=0}^{\binom{n}{2}} i_{B'}(n,k) = (2n-1)!!, \end{equation*} 
where $(2n-1)!!=1\times 3 \times\cdots \times (2n-1)$.
\end{corollary}
Consider the set of $2n-1$ symbols 
$$\Sigma'_{2,n}:=\left\{1,2,\overline{2},\ldots,n,\overline{n}\right\}.$$

Let $C_j=\{\sigma \in B'_n: \sigma(1)=j\}$ for $j\in \Sigma'_{2,n}$. It is clear that
\begin{equation}\label{tmr}
|C_j|=(2n-3)!!, 
\end{equation}
with $|C_1|=1$ when $n=1$.

\medskip

The set $C_j$ gives the decomposition
\[B'_n=\biguplus_{j\in \Sigma'_{2,n}}C_j.\]
Therefore, we have
\begin{equation}\label{bnn}
\mathcal{B}'_n=\sum_{\omega\in B'_n}inv(\omega)=\sum_{j\in \Sigma'_{2,n}}\sum_{\sigma \in C_j}inv(\sigma),
\end{equation}
where $\mathcal{B}'_n$ represents the total number of inversions of all overlined permutations in $B'_n$. 

\medskip

Let $\pi=\left(
  \begin{array}{cccc}
    1 &2&\cdots & n \\
    j &\pi(2) &\cdots & \pi(n)\\
  \end{array}\right)\in C_j$ and $\tau$ be is an overlined permutation in $B'_{n-1}$ defined by 
\[
\tau=\left(
  \begin{array}{ccc}
    a_1 & \cdots & a_{n-1} \\
    \pi(2) & \cdots & \pi(n) \\
  \end{array}\right)\in B'\left([n]\backslash\{|j|\}\right),
\]
where $a_1,\ldots,a_{n-1}$ are an arrangement of elements of $[n]\backslash\{|j|\}$ in increasing order, and $B'\left([n]\backslash\{|j|\}\right)$ is the group of all the overlined permutation of the set $[n]\backslash\{|j|\}$. So, if we set 

\[
\pi_{\tau,j}=\left(
  \begin{array}{cccc}
    1 &2  &\cdots& n\\
    j &\tau(a_1) & \cdots & \tau(a_{n-1}) \\
   \end{array} \right),
\]
then we obtain $\pi=\pi_{\tau,j}$. Hence by the definition of the statistic of the classical inversion, we conclude 
\begin{equation}\label{eqt}
inv(\pi)=(j-1)+inv(\tau).
\end{equation}
Equation \eqref{eqt} give us a recursive formula for $\mathcal{B}'_n$ and we state it as a proposition.
\begin{proposition}
We have $\mathcal{B}'_1=0$, and 
\begin{equation}\label{eip}
\mathcal{B}'_n=(2n-3)!!n(n-1)+(2n-1)\mathcal{B}'_{n-1}, \text{\ for\ }n\geq 2.
\end{equation} 
\end{proposition}
\begin{proof}
Since $B'_1=\{\iota\}$, then $\mathcal{B}'_1=inv(\iota)=0$. Now suppose $n\geq 2$. We have two cases:

\medskip

\noindent{\bf Case 1.} For $j=1$, we obtain
\begin{align*}
\sum_{\pi\in C_j}inv(\pi)&=\sum_{\tau\in B'\left([n]\backslash\{|j|\}\right)}inv(\pi_{\tau,j})\\
&=\sum_{\tau\in B'\left([n]\backslash\{|j|\}\right)}inv(\tau)\\
&=\mathcal{B}'_{n-1}.
\end{align*}
\noindent{\bf Case 2.} For $j\in \Sigma'_{2,n}\backslash \{1\}$. From equation \eqref{tmr}, we obtain
\begin{align*}
\sum_{\pi\in C_j}inv(\pi)&=\sum_{\tau\in B'\left([n]\backslash\{|j|\}\right)}inv(\pi_{\tau,j})\\
&=\sum_{\tau\in B'\left([n]\backslash\{|j|\}\right)}[(j-1)+inv(\tau)]\\
&=(j-1)(2n-3)!!+\mathcal{B}'_{n-1}.
\end{align*}
From equation \eqref{bnn}, we get
\begin{align*}
\mathcal{B}'_{n}&=\mathcal{B}'_{n-1}+2\sum_{j=2}^{n}\left[(j-1)(2n-3)!!+\mathcal{B}'_{n-1}\right]\\
&=2(2n-3)!!\sum_{j=2}^{n}(j-1)+(2n-1)\mathcal{B}'_{n-1}\\
&=(2n-3)!!n(n-1)+(2n-1)\mathcal{B}'_{n-1}
\end{align*}
as desired.
\end{proof}

Note that $\mathcal{B}'_{n}=\sum_{k=0}^{\binom{n}{2}}inv_B'(n,k)k$. The numbers $i_{B'}(n,k)$ form a triangle known as the "\textbf{over-Mahonian triangle}", as shown in Table 1. 

\medskip

\begin{center}
$
\begin{array}{|c|c|ccccccccccc|}
\hline
n\backslash k & \mathcal{B}'_{n} &\ \ 0 \ \  & \ \ 1 \ \  & \ \ 2 \ \  & \ \ 3 \ \  & \ \ 4 \
\  & \ \ 5 \ \  & \ \ 6 \ \  & \ \ 7 \ \  & \ \ 8 \ \  & \ \ 9 \ \  & \ \ 10
      \\ 
\hline
1 & 0 & 1 &  &  &  &  &  &  &  &  &  &      \\ 
2 & 2 & 1 & 2 &  &  &  &  &  &  &  &  &        \\ 
3 & 28 & 1 & 4 & 6 & 4 &  &  &  &  &  &  &        \\ 
4 & 376 & 1 & 6 & 16 & 26 & 28 & 20 & 8 &  &  &  &       \\ 
5 & & 1 & 8 & 30 & 72 & 126 & 172 & 188 & 164 & 112 & 56 & 16      \\
\hline
\end{array}
$
\end{center}

\begin{center}
Table 1: The over-Mahonian triangle. 
\end{center}

\section{Combinatorial interpretations of the over-Mahonian numbers}\label{sec:3}
In this section, we provide combinatorial interpretations of the over-Mahonian numbers through lattice paths, overpartitions, and tilings.
\subsection{A lattice path interpretation}\label{path}
Ghemit and Ahmia in \cite{GA1} demonstrated that the Mahonian number $i(n,k)$ counts the number of lattice paths from $ u_{1}=(0,0) $ to $ v_{1}=(n-1,k) $ that take at most $j$ North steps at the level $j$.

\medskip

Following a similar approach, we provide a combinatorial interpretation for the over-Mahonian numbers $i_{B'}(n,k)$ using lattice paths.

\medskip

\begin{figure}[ht]

\begin{center}
\begin{tikzpicture}

\draw[step=0.8cm,color=black!30] (-0.81,0) grid (1.6,4.8);
\draw [line width=1pt](-0.8,0) -- (0,0) -- (0,0.8) -- (0.8,1.6) -- (0.8,2.4)--(1.6,3.2)--(1.6,3.2)--(1.6,4)--(1.6,4.8);

\fill[blue] (-0.8,0) circle(1.9pt) ;
\fill[blue] (1.6,4.8) circle(1.9pt) ;
\node[right=0.1pt] at (1.7,4.8){$v_1$};
\node at (-1.1,0){$u_1$};
\node at (-1.8,-0.5){$levels$};
\node at (-0.8,-0.5){$0$};
\node at (0,-0.5){$1$};
\node at (0.8,-0.5){$2$};
\node at (1.6,-0.5){$3$};

\end{tikzpicture}
\end{center}

\caption{A path $P$ in $\mathcal{P}^{B' }_{4,6}$.}
  \label{fg1}
\end{figure}
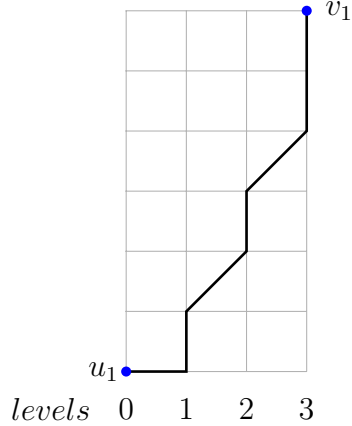

Let $\mathcal{P}^{B'}_{n,k}$ represent the set of lattice paths from the point $(0,0)$ to $(n-1,k)$, using only North steps (vertical steps $(0,1)$), East steps (horizontal steps $(1,0)$), and North-East steps (diagonal steps $(1,1)$). These paths must satisfy the condition that, for each vertical level $j$, the number of North steps is at most $j$ unless a diagonal step occurs before level $j$, in which case the number of North steps is restricted to at most $j-1$ after the diagonal step. The levels correspond to the vertical lines ranging from $0$ to $n-1$, as illustrated in Figure \ref{fg1}.

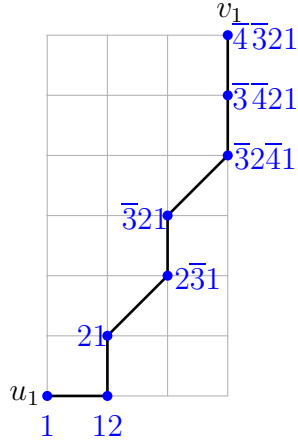
\begin{figure}[ht]

	\begin{center}
		\begin{tikzpicture}

		\draw[step=0.8cm,color=black!30] (-0.81,0) grid (1.6,4.8);
		\draw [line width=1pt](-0.8,0) -- (0,0) -- (0,0.8) -- (0.8,1.6) -- (0.8,2.4)--(1.6,3.2)--(1.6,3.2)--(1.6,4)--(1.6,4.8);
		
		\fill[blue] (-0.8,0) circle(1.9pt) ;
		\fill[blue] (0,0) circle(1.9pt) ;
		\fill[blue] (0,0.8) circle(1.9pt) ;
		\fill[blue] (0.8,2.4) circle(1.9pt) ;
		\fill[blue] (1.6,3.2) circle(1.9pt) ;
		\fill[blue] (0.8,1.6) circle(1.9pt) ;
		\fill[blue] (1.6,4) circle(1.9pt) ;
		\fill[blue] (1.6,4.8) circle(1.9pt) ;
		\node[right=0.1pt] at (1.3,5.1){$v_1$};
		\node at (-1.1,0){$u_1$};

		\node[blue,rectangle] at (-0.8,-0.4) {$1$};
		\node[blue,rectangle] at (0.0,-0.4) {$12$};
	\node[blue,rectangle] at (-0.2,0.8) {$21$};
	\node[blue,rectangle] at (1.2,1.6) {$2\overline{3}1$};
	\node[blue,rectangle] at (0.5,2.4) {$\overline{3}21$};
	
	\node[blue,rectangle] at (2.1,3.2) {$\overline{3}2\overline{4}1$};
	\node[blue,rectangle] at (2.1,4) {$\overline{3}\hspace{0.05cm} \overline{4}21$};
	\node[blue,rectangle] at (2.1,4.8) {$\overline{4}\hspace{0.05cm}\overline{3}21$};
		\end{tikzpicture}
	\end{center}
	

	\caption{ The path associated to the permutation $\sigma =\overline{4}\hspace{0.05cm} \overline{3}21$.}
	\label{fg2}
\end{figure}

\begin{theorem}
The over-Mahonian number $i_{B'}(n,k)$ represents the count of lattice paths from $ u_{1}=(0,0) $ to $ v_{1}=(n-1,k) $, where the path takes at most $ j $ North steps at level $ j $ before any diagonal step, and at most $(j-1)$ North steps at level $j$ after a diagonal step. Specifically, we have
$$i_{B'}(n,k)= \mid \mathcal{P}^{B'}_{n,k} \mid.$$
 \end{theorem}
\begin{proof}
Since $i_{B'}(n,k)$ counts the number of permutations of length $n$ with exactly $k$ overlined  inversions, it suffices to show a bijection between these permutations and the lattice paths in $\mathcal{P}^{B'}_{n,k}$.  We do so as follows:

\smallskip

\noindent For each path $P\in \mathcal{P}^{B'}_{n,k}$, we can uniquely determine the associated permutation. This is done by assigning the entry $1$ to the point $(0,0)$, and then following the steps of the path. The first step must be either an East step or a North-East step. If the step is East, we move to the point $(1,0)$ and place the entry $2$ to the right of $1$ (resulting in the partial permutation $12$). If the step is North-East, we move to $(1,1)$ and place the entry $ \overline{2} $  to the left of $ 1 $, yielding the partial permutation $ \overline{2}1$, which introduces an inversion. 

\medskip

\noindent \textbf{At the point $(1,0)$}, there are three possible cases for the next step of $P$: 

\begin{description}
\item [$\bullet$] If the next step is East, we add the entry $3$ to the right of $2$, resulting in the partial permutation $123$.
\item[$\bullet$] If the next step is North-East, we add the entry $\overline{3}$ to the right of $ 2 $ and shift $\overline{3}$ one position to the left, producing the permutation $ 1\overline{3}2 $, which introduces another inversion.
\item[$\bullet$] If the next step is North, we shift the entry $2$ one position to the left, yielding the permutation $21$, which also introduces an inversion.
\end{description}

\noindent \textbf{At the point $(1,1)$}, there are two possible cases for the next step: 
\begin{description}
\item[$\bullet$] If the next step is East, we add the entry $3$ to the right of $1$, producing the permutation $ \overline{2}13$.
\item[$\bullet$] If the next step is North-East, we add $\overline{3}$  to the right of $1$ and shift $\overline{3}$ one position to the left, resulting in the permutation $ \overline{2}\hspace{0.05cm}\overline{3}1 $ which introduces an inversion.
\end{description}
We continue applying these rules at each step of the path until we reach the point $(n-1,k)$, which gives us the desired permutation. An example of this process is shown in Figure \ref{fg2}.\end{proof}

\subsection{An overpartition interpretation}\label{overp}
This subsection explores combinatorial interpretations of the over-Mahonian numbers using overpartitions.

\medskip

We begin with two fundamental definitions:
\begin{definition} A \textbf{partition} $\lambda = (\lambda_1, \lambda_2, \ldots, \lambda_k)$ of a number $n$ is defined as a non-increasing sequence of positive integers, i.e., $\lambda_1 \geq \lambda_2 \geq \cdots \geq \lambda_k$, where the sum of the sequence equals $n$. The function $p(n)$ represents the number of such partitions of $n$, with the convention that $p(0) = 1$. 
\end{definition}

\noindent For $1 \leq i \leq k$, $\lambda_i$ is referred to as a part of $\lambda$. While partitions typically consist of positive integers, we may sometimes allow "zero" as a part for specific purposes.

\medskip

\noindent The \textit{length} of $\lambda$, denoted by $l(\lambda)$, is the number of parts in $\lambda$, and the \textit{weight} of $\lambda$, denoted by $|\lambda|$, is the total sum of its parts.

\medskip

\begin{definition}\cite[Corteel and Lovejoy]{Cort} An \textbf{overpartition} of $n$ is a non-increasing sequence of natural numbers whose sum is $n$, where the first occurrence (or the last occurrence) of any number may be overlined. The function $\overline{p}(n)$ represents the number of overpartitions of $n$, with the convention that $\overline{p}(0) = 1$. 
\end{definition}

For example, there are $8$ overpartitions of $3$, enumerated as follows: $$3,\overline{3},2+1,\overline{2}+1,2+\overline{1},\overline{2}+\overline{1},1+1+1,\overline{1}+1+1.$$

\medskip

Since the bijection between lattice paths and tilings preserves weight, we derive the following interpretation of the over-Mahonian numbers in terms of overpartitions:
\begin{theorem} 
The over-Mahonian number $i_{B'}(n,k)$  represents the count of overpartitions into $k$ parts, where each part $j$ can appear at most $j$ times if it is not overlined, and at most $j-1$ times if it is overlined. Additionally, the largest part must satisfy $\leq n-1$. 
\end{theorem}
\begin{proof}
The term $i_{B'}(n,k)$ represents the count of lattice paths from $(0,0)$ to $(n-1,k)$, where the paths adhere to specific rules: at most $j$ North steps are allowed at level $j$ if no diagonal step occurs beforehand, and at most $j-1$ North steps are permitted at level $j$ if a diagonal step has occurred. If we interpret the parts (resp. the overlined parts) as representing the cases above each North step (resp. North-East step) in the paths associated with $i_{B'}(n,k)$, we can deduce a direct connection to overpartitions. Specifically, $i_{B'}(n,k)$ enumerates the number of overpartitions of $k$ into parts where each part $j$ can appear at most $j$ times if it is not overlined and at most $j-1$ times if it is overlined, with the additional condition that the largest part must be at most $n-1$.
\end{proof}

\begin{figure}[ht!]	
	\begin{center}
		\begin{tikzpicture}
		\draw[step=0.8cm,color=black!30] (-0.81,0) grid (1.6,1.6);
		\draw [line width=1pt](-0.8,0) -- (0,0) -- (0,0.8) -- (0.8,0.8) -- (0.8,1.6)--(1.6,1.6);
		\fill[blue] (-0.8,0) circle(1.9pt) ;
		\fill[blue] (1.6,1.6) circle(1.9pt) ;
		\node[right=0.1pt] at (1.7,1.6){$v_1$};
		\node at (-1.1,0){$u_1$};
		\node at (0.5,-0.5){$\lambda=(1,2)$};
		\end{tikzpicture}
		\begin{tikzpicture}
		\draw[step=0.8cm,color=black!30] (-0.81,0) grid (1.6,1.6);
		\draw [line width=1pt](-0.8,0) -- (0,0) -- (0,0.8) -- (0.8,0.8) -- (1.6,0.8)--(1.6,1.6);
		\fill[blue] (-0.8,0) circle(1.9pt) ;
		\fill[blue] (1.6,1.6) circle(1.9pt) ;
		\node[right=0.1pt] at (1.7,1.6){$v_1$};
		\node at (-1.1,0){$u_1$};
		\node at (0.5,-0.5){$\lambda=(1,3)$};
		\end{tikzpicture}
		\begin{tikzpicture}
		\draw[step=0.8cm,color=black!30] (-0.81,0) grid (1.6,1.6);
		\draw [line width=1pt](-0.8,0) -- (0,0) -- (0.8,0) -- (0.8,0.8) -- (1.6,0.8)--(1.6,1.6);
		
		\fill[blue] (-0.8,0) circle(1.9pt) ;
		\fill[blue] (1.6,1.6) circle(1.9pt) ;
		\node[right=0.1pt] at (1.7,1.6){$v_1$};
		\node at (-1.1,0){$u_1$};
		\node at (0.5,-0.5){$\lambda=(2,3)$};
	
		\end{tikzpicture}
		\begin{tikzpicture}
		\draw[step=0.8cm,color=black!30] (-0.81,0) grid (1.6,1.6);
		\draw [line width=1pt](-0.8,0) -- (0,0) -- (0.8,0) -- (0.8,0.8) -- (0.8,1.6)--(1.6,1.6);
		
		\fill[blue] (-0.8,0) circle(1.9pt) ;
		\fill[blue] (1.6,1.6) circle(1.9pt) ;
		\node[right=0.1pt] at (1.7,1.6){$v_1$};
		\node at (-1.1,0){$u_1$};
		\node at (0.5,-0.5){$\lambda=(2,2)$};
	
		\end{tikzpicture}
		\begin{tikzpicture}
		\draw[step=0.8cm,color=black!30] (-0.81,0) grid (1.6,1.6);
		\draw [line width=1pt](-0.8,0) -- (0,0.8) -- (0.8,0.8) -- (0.8,1.6)--(1.6,1.6);
		\fill[blue] (-0.8,0) circle(1.9pt) ;
		\fill[blue] (1.6,1.6) circle(1.9pt) ;
		\node[right=0.1pt] at (1.7,1.6){$v_1$};
		\node at (-1.1,0){$u_1$};
		\node at (0.5,-0.5){$\lambda=(\overline{1},2)$};
		\end{tikzpicture}
		\begin{tikzpicture}
		\draw[step=0.8cm,color=black!30] (-0.81,0) grid (1.6,1.6);
		\draw [line width=1pt](-0.8,0) -- (0,0.8) -- (0.8,0.8) -- (1.6,0.8)--(1.6,1.6);
		\fill[blue] (-0.8,0) circle(1.9pt) ;
		\fill[blue] (1.6,1.6) circle(1.9pt) ;
		\node[right=0.1pt] at (1.7,1.6){$v_1$};
		\node at (-1.1,0){$u_1$};
		\node at (0.5,-0.5){$\lambda=(\overline{1},3)$};
		\end{tikzpicture}
		\begin{tikzpicture}
		\draw[step=0.8cm,color=black!30] (-0.81,0) grid (1.6,1.6);
		\draw [line width=1pt](-0.8,0) -- (0,0) -- (0.8,0.8) -- (1.6,0.8)--(1.6,1.6);
		
		\fill[blue] (-0.8,0) circle(1.9pt) ;
		\fill[blue] (1.6,1.6) circle(1.9pt) ;
		\node[right=0.1pt] at (1.7,1.6){$v_1$};
		\node at (-1.1,0){$u_1$};
		\node at (0.5,-0.5){$\lambda=(\overline{2},3)$};
	
		\end{tikzpicture}
		\begin{tikzpicture}
		\draw[step=0.8cm,color=black!30] (-0.81,0) grid (1.6,1.6);
		\draw [line width=1pt](-0.8,0) -- (0,0)  -- (0.8,0.8) -- (0.8,1.6)--(1.6,1.6);
		
		\fill[blue] (-0.8,0) circle(1.9pt) ;
		\fill[blue] (1.6,1.6) circle(1.9pt) ;
		\node[right=0.1pt] at (1.7,1.6){$v_1$};
		\node at (-1.1,0){$u_1$};
		\node at (0.5,-0.5){$\lambda=(\overline{2},2)$};
	
		\end{tikzpicture}
		\begin{tikzpicture}
		\draw[step=0.8cm,color=black!30] (-0.81,0) grid (1.6,1.6);
		\draw [line width=1pt](-0.8,0) -- (0,0) -- (0,0.8) --  (0.8,1.6)--(1.6,1.6);
		\fill[blue] (-0.8,0) circle(1.9pt) ;
		\fill[blue] (1.6,1.6) circle(1.9pt) ;
		\node[right=0.1pt] at (1.7,1.6){$v_1$};
		\node at (-1.1,0){$u_1$};
		\node at (0.5,-0.5){$\lambda=(1,\overline{2})$};
		\end{tikzpicture}
		\begin{tikzpicture}
		\draw[step=0.8cm,color=black!30] (-0.81,0) grid (1.6,1.6);
		\draw [line width=1pt](-0.8,0) -- (0,0) -- (0,0.8) -- (0.8,0.8) -- (1.6,1.6);
		\fill[blue] (-0.8,0) circle(1.9pt) ;
		\fill[blue] (1.6,1.6) circle(1.9pt) ;
		\node[right=0.1pt] at (1.7,1.6){$v_1$};
		\node at (-1.1,0){$u_1$};
		\node at (0.5,-0.5){$\lambda=(1,\overline{3})$};
		\end{tikzpicture}
		\begin{tikzpicture}
		\draw[step=0.8cm,color=black!30] (-0.81,0) grid (1.6,1.6);
		\draw [line width=1pt](-0.8,0) -- (0,0) -- (0.8,0) -- (0.8,0.8) --(1.6,1.6);
		
		\fill[blue] (-0.8,0) circle(1.9pt) ;
		\fill[blue] (1.6,1.6) circle(1.9pt) ;
		\node[right=0.1pt] at (1.7,1.6){$v_1$};
		\node at (-1.1,0){$u_1$};
		\node at (0.5,-0.5){$\lambda=(2,\overline{3})$};
	
		\end{tikzpicture}
		\begin{tikzpicture}
		\draw[step=0.8cm,color=black!30] (-0.81,0) grid (1.6,1.6);
		\draw [line width=1pt](-0.8,0) -- (0,0) -- (0.8,0) -- (1.6,0) -- (1.6,0.8)--(1.6,1.6);
		
		\fill[blue] (-0.8,0) circle(1.9pt) ;
		\fill[blue] (1.6,1.6) circle(1.9pt) ;
		\node[right=0.1pt] at (1.7,1.6){$v_1$};
		\node at (-1.1,0){$u_1$};
		\node at (0.5,-0.5){$\lambda=(3,3)$};
	
		\end{tikzpicture}
		\begin{tikzpicture}
		\draw[step=0.8cm,color=black!30] (-0.81,0) grid (1.6,1.6);
		\draw [line width=1pt](-0.8,0) -- (0,0.8) -- (0.8,1.6)--(1.6,1.6);
		\fill[blue] (-0.8,0) circle(1.9pt) ;
		\fill[blue] (1.6,1.6) circle(1.9pt) ;
		\node[right=0.1pt] at (1.7,1.6){$v_1$};
		\node at (-1.1,0){$u_1$};
		\node at (0.5,-0.5){$\lambda=(\overline{1},\overline{2})$};
		\end{tikzpicture}
		\begin{tikzpicture}
		\draw[step=0.8cm,color=black!30] (-0.81,0) grid (1.6,1.6);
		\draw [line width=1pt](-0.8,0) -- (0,0.8) -- (0.8,0.8)--(1.6,1.6);
		\fill[blue] (-0.8,0) circle(1.9pt) ;
		\fill[blue] (1.6,1.6) circle(1.9pt) ;
		\node[right=0.1pt] at (1.7,1.6){$v_1$};
		\node at (-1.1,0){$u_1$};
		\node at (0.5,-0.5){$\lambda=(\overline{1},\overline{3})$};
		\end{tikzpicture}
		\begin{tikzpicture}
		\draw[step=0.8cm,color=black!30] (-0.81,0) grid (1.6,1.6);
		\draw [line width=1pt](-0.8,0) -- (0,0)-- (0.8,0.8)--(1.6,1.6);
		
		\fill[blue] (-0.8,0) circle(1.9pt) ;
		\fill[blue] (1.6,1.6) circle(1.9pt) ;
		\node[right=0.1pt] at (1.7,1.6){$v_1$};
		\node at (-1.1,0){$u_1$};
		\node at (0.5,-0.5){$\lambda=(\overline{2},\overline{3})$};
	
		\end{tikzpicture}
		\begin{tikzpicture}
		\draw[step=0.8cm,color=black!30] (-0.81,0) grid (1.6,1.6);
		\draw [line width=1pt](-0.8,0) -- (0,0) -- (0.8,0) -- (1.6,0.8)--(1.6,1.6);
		
		\fill[blue] (-0.8,0) circle(1.9pt) ;
		\fill[blue] (1.6,1.6) circle(1.9pt) ;
		\node[right=0.1pt] at (1.7,1.6){$v_1$};
		\node at (-1.1,0){$u_1$};
		\node at (0.5,-0.5){$\lambda=(\overline{3},3)$};
	
		\end{tikzpicture}

	\end{center}



\caption{All lattice paths/overpartitions corresponding to $n=4$ and $ k=2$.}
  \label{ff2}
\end{figure}
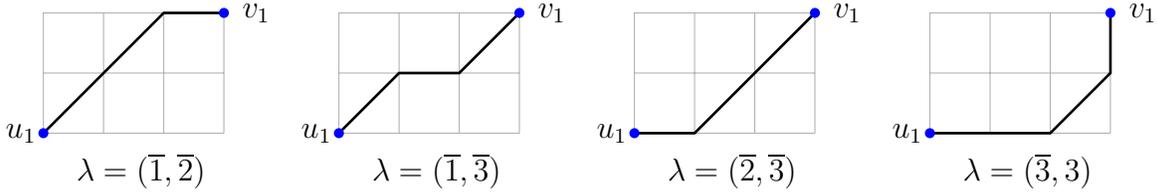

For instance, Figure \ref{ff2} illustrates the $16$ overpartitions with $k=2$ parts, where the largest part is $n-1=3$.

\subsection{A tiling interpretation} \label{tile}
Ghemit and Ahmia \cite{GA2} demonstrated that the Mahonian number $i(n,k)$ represents the number of ways to tile a board of length $n+k-1$ using $k$ red squares and $n-1$ blue squares, subject to the condition that no more than $j$ consecutive red squares appear if preceded by $j$ blue squares. Building on their findings, this subsection provides a combinatorial interpretation of the over-Mahonian numbers using the same combinatorial structures.

\medskip

\begin{figure}[ht]

\begin{center}
\begin{tikzpicture}

\draw[step=0.8cm,color=black!30] (0,0) grid(3.2,0.8);

\fill[fill=blue!60,draw=black!50] (0.05,0.05) rectangle (0.75,0.75);
\fill[fill=red!90,draw=black!50] (0.85,0.05) rectangle (1.55,0.75);
\fill[fill=black!90,draw=black!50] (1.65,0.05) rectangle (2.35,0.75);
\fill[fill=black!90,draw=black!50] (2.35,0.05) rectangle (3.15,0.75);

\end{tikzpicture}
\hspace{1.2cm}
\begin{tikzpicture}

\draw[step=0.8cm,color=black!30] (-0.81,0) grid (0.8,1.6);
\draw [line width=1pt](-0.8,0) -- (0,0) -- (0,0.8) -- (0.8,1.6) ;

\fill[blue] (-0.8,0) circle(1.9pt) ;
\fill[blue] (0.8,1.6) circle(1.9pt) ;
\node[right=0.1pt] at (0.9,1.6){$v_1$};
\node at (-1.1,0){$u_1$};

	\node[blue,rectangle] at (-0.8,-0.4) {$1$};
	\node[blue,rectangle] at (0.0,-0.4) {$12$};
	\node[blue,rectangle] at (-0.2,0.8) {$21$};
	\node[blue,rectangle] at (0.7,1.9) {$2\overline{3}1$};
\end{tikzpicture}.
\end{center}
  \caption{The tiling/lattice path corresponding to the permutation $\sigma= \hspace{0.05cm} 2\overline{3}1$.}\label{tqq}
\end{figure}
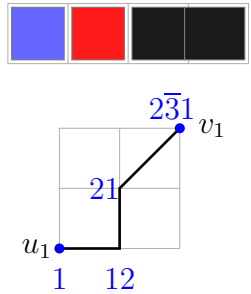

Let $\mathcal{T}^{B'}_{n,k} $ denote the set of all tilings of a board of length $(n+k-1)$ that satisfy the following conditions: 
\begin{enumerate}
\item The tiling uses at most $n-1$ blue squares, at most $k$ red squares, and black rectangles (each covering an area of two squares, referred to as square pairs).
\item The number of consecutive red squares does not exceed $j$ if preceded by $j$ blue squares, and does not exceed $j-1$ if preceded by a black rectangle.
\item The total number of red squares and black rectangles is exactly $k$.
\item The total number of blue squares and black rectangles is exactly $n-1$.
\item Each black rectangle is treated as equivalent to one blue square followed by one red square in that order.
\end{enumerate}
\begin{remark}
There exists a bijection between the sets $\mathcal{P}^{B'}_{n,k} $ and $\mathcal{T}^{B'}_{n,k} $. Specifically, each East step corresponds to a blue square, each North step corresponds to a red square, and each diagonal step corresponds to a black rectangle. This correspondence is illustrated in Figure \ref{tqq}.
\end{remark}

Based on this observation, we can derive the following tiling interpretation for the over-Mahonian numbers:
\begin{corollary} The over-Mahonian number $i_{B'}(n,k)$ is precisely the cardinality of the set $\mathcal{T}^{B'}_{n,k}$. In other words, 
     $i_{B'}(n,k)=\mid\mathcal{T}^{B'}_{n,k}\mid. $ 
\end{corollary}
\begin{example}
For $n=3$ and $k=2$, the $i_{B'}(3,2)$  tilings consist of $2$ blue squares, $2$ red squares, and a black rectangle, as illustrated in Figure \ref{ti}.

\begin{figure}[ht]

\begin{center}
\begin{tikzpicture}

\draw[step=0.8cm,color=black!30] (0,0) grid(3.2,0.8);

\fill[fill=blue!90,draw=black!50] (0.05,0.05) rectangle (0.75,0.75);
\fill[fill=blue!90,draw=black!50] (0.85,0.05) rectangle (1.55,0.75);
\fill[fill=red!90,draw=black!50] (1.65,0.05) rectangle (2.30,0.75);
\fill[fill=red!90,draw=black!50] (2.40,0.05) rectangle (3.15,0.75);

\end{tikzpicture}
\hspace{2cm}
\begin{tikzpicture}

\draw[step=0.8cm,color=black!30] (0,0) grid(3.2,0.8);

\fill[fill=blue!90,draw=black!50] (0.05,0.05) rectangle (0.75,0.75);
\fill[fill=red!90,draw=black!50] (0.85,0.05) rectangle (1.55,0.75);
\fill[fill=blue!90,draw=black!50] (1.65,0.05) rectangle (2.30,0.75);
\fill[fill=red!90,draw=black!50] (2.40,0.05) rectangle (3.15,0.75);

\end{tikzpicture}
\hspace{2cm}
\begin{tikzpicture}

\draw[step=0.8cm,color=black!30] (0,0) grid(3.2,0.8);

\fill[fill=blue!90,draw=black!50] (0.05,0.05) rectangle (0.75,0.75);
\fill[fill=black!90,draw=black!50] (0.85,0.05) rectangle (1.65,0.75);
\fill[fill=black!90,draw=black!50] (1.65,0.05) rectangle (2.30,0.75);
\fill[fill=red!90,draw=black!50] (2.40,0.05) rectangle (3.15,0.75);

\end{tikzpicture} 
\begin{tikzpicture}

\draw[step=0.8cm,color=black!30] (0,0) grid(3.2,0.8);

\fill[fill=blue!90,draw=black!50] (0.05,0.05) rectangle (0.75,0.75);
\fill[fill=red!90,draw=black!50] (0.85,0.05) rectangle (1.55,0.75);
\fill[fill=black!90,draw=black!50] (1.65,0.05) rectangle (2.35,0.75);
\fill[fill=black!90,draw=black!50] (2.35,0.05) rectangle (3.15,0.75);

\end{tikzpicture}
\hspace{2cm}
\begin{tikzpicture}

\draw[step=0.8cm,color=black!30] (0,0) grid(3.2,0.8);

\fill[fill=black!90,draw=black!50] (0.05,0.05) rectangle (0.85,0.75);
\fill[fill=black!90,draw=black!50] (0.85,0.05) rectangle (1.55,0.75);
\fill[fill=blue!90,draw=black!50] (1.65,0.05) rectangle (2.30,0.75);
\fill[fill=red!90,draw=black!50] (2.40,0.05) rectangle (3.15,0.75);

\end{tikzpicture}
\hspace{2cm}
\begin{tikzpicture}

\draw[step=0.8cm,color=black!30] (0,0) grid(3.2,0.8);

\fill[fill=black!90,draw=black!50] (0.05,0.05) rectangle (0.85,0.75);
\fill[fill=black!90,draw=black!50] (0.85,0.05) rectangle (1.55,0.75);
\fill[fill=black!90,draw=black!50] (1.65,0.05) rectangle (2.35,0.75);
\fill[fill=black!90,draw=black!50] (2.35,0.05) rectangle (3.15,0.75);

\end{tikzpicture}

\end{center}
  \caption{The tiling interpretation associated to  $i_{B'}(3,2)$.}\label{ti}
\end{figure}
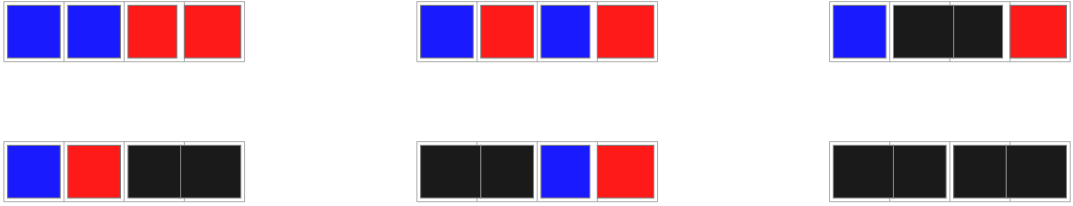
\end{example}


\section{Combinatorial proof of the log-concavity and unimodality of the over-Mahonian numbers}\label{sec:4}
Log-concave and unimodal sequences frequently appear in combinatorics, geometry, and algebra. For a comprehensive overview of the various techniques used to study sequences and polynomials that are log-concave or unimodal, the reader is referred to \cite{bren1, stan}.

\medskip

A sequence of nonnegative numbers $\lbrace x_{k}\rbrace_{k}$ is said to be log-concave if it satisfies the inequality $x_{i-1}x_{i+1}\leq x_{i}^{2}$ for all $i>0$,
. This condition is equivalent to the following for relevant results (see \cite{bren2, stan}):
\begin{equation*}
x_{i-1}x_{j+1}\leq x_{i}x_{j}\hspace{0.5cm} for \hspace{0.5cm} j\geq i\geq1.
\end{equation*}

\medskip

A finite sequence of real numbers $a_{0},\ldots,a_{m}$ is called unimodal if there exists an index $0\leq m^{*}\leq m $, known as the mode of the sequence, such that the sequence first increases up to $k=m^{*}$ and then decreases thereafter. Specifically, the sequence satisfies  $a_{0}\leq a_{1}\leq \cdots\leq a_{m^{*}} $ and $a_{m^{*}}\geq a_{m^{*}+1}\geq \cdots \geq a_{m}$. It is clear that any log-concave sequence is unimodal \cite{bren2}.
 
\medskip 
 
A polynomial is called log-concave (resp. unimodal) if its coefficients form a log-concave (resp. unimodal) sequence. Notably, Gaussian polynomials and $q$-multinomial coefficients are well-known examples of unimodal sequences, as established in \cite[Theorem 3.11]{Syl}. Additionally, it is a classic result, with proofs available in works such as \cite{bona}, that the product of log-concave (resp. unimodal) polynomials remains log-concave (resp. unimodal). For instance, the polynomial $\sum_{k=0}^{\binom{n}{2}}i(n,k)x^k$ is log-concave (resp. unimodal), meaning the sequence $i(n,0), \ldots, i(n,\binom{n}{2})$ is log-concave (resp. unimodal).

\medskip

B\'ona \cite{bonalog} provided the first non-generating function proof of the log-concavity of the Mahonian numbers $i(n,k)$ in $k$, using an injection property and an induction hypothesis over  n, although the injection was non-constructive. To address this, Ghemit and Ahmia \cite{GA2} offered a constructive proof that replaces B\'ona's non-constructive injection.

\medskip

The proof of the log-concavity (or unimodality) of the over-Mahonian numbers $\{i_{B'}(n,k)\}_{k}$ in $k$ can be easily derived from Theorem \ref{thm2}, since the generating function $\sum_{k=0}^{\binom{n}{2}}i_{B'}(n,k) x^k$ is the product of log-concave (or unimodal) polynomials. However, in this section, motivated by the work of Ghemit and Ahmia \cite{GA2}, we provide a new approach to demonstrate that the sequence of over-Mahonian numbers $\{i_{B'}(n,k)\}_{k}$ s log-concave in $k$, and therefore unimodal, by constructing an appropriate injection.

\bigskip

We begin by introducing the following key definition.
\begin{definition}
Let $\sigma$ be a permutation of length $n$. For $ 0\leq i\leq n-1$, let $m_{i}^{\sigma}$  denote the number of times the entry  $i + 1$ appears as the first element of an overlined inversion of $\sigma$. The total number of overlined inversions of  $\sigma$ is then given by $\sum_{i=0}^{n-1}m_{i}^{\sigma}$.
\end{definition}

Based on this definition, we can derive the following result regarding the log-concavity of over-Mahonian numbers:
\begin{theorem}
The sequence of over-Mahonian numbers $\{i_{B'}(n,k)\}_{k}$ is log-concave in $k$. Specifically, we have the inequality: 
\begin{equation*}
\left( i_{B'}(n,k)\right) ^{2}-i_{B'}(n,k-1)i_{B'}(n,k+1)\geq0,
\end{equation*} 
for all $0<k\leq\binom{n}{2}$.
\end{theorem}

\begin{proof}
Let $I_{B'}(n,k)$ denote the set of permutations of length $n$ with $k$ overlined inversions, and let $i_{B'}(n,k)=|I_{B'}(n,k)|$ represent the number of such permutations.

\medskip

To prove the log-concavity of the sequence of over-Mahonian numbers ${_{B'}(n,k)}_{k}$, it is equivalent to find an injection $f_{n,k,k}$ from $I_{B'}(n,k+1)\times I_{B'}(n,k-1)$  to $I_{B'}(n,k)\times I_{B'}(n,k)$, where $f_{n,k,k}(\sigma,\tau)=(\theta,\pi)$ with $(\sigma,\tau)\in I_{B'}(n,k+1)\times I_{B'}(n,k-1)$ and $(\theta,\pi)\in I_{B'}(n,k)\times I_{B'}(n,k)$.

\bigskip

Let $(\sigma,\tau) \in I_{B'}(n,k+1) \times I_{B'}(n,k-1)$.
We will define $I$ as the largest integer which satisfies the following three conditions:
\medskip
\begin{align}
	&\sum_{j\geq I}m^\sigma_{j}\geq \sum_{j\geq I+1}m^\tau_{j}+1,\text{\quad if the entry }(j+1)=(I+1) \text{\quad of the  }\sigma\text{\quad not overlined},
\end{align}\label{cond1}
\begin{align}
	&\sum_{j\geq I}m^\sigma_{j}> \sum_{j\geq I+1}m^\tau_{j}+1,\text{\quad if the entry }(j+1)=(I+1) \text{\quad of the  } \sigma\text{\quad is overlined}\label{cond2}
\end{align}
and
\begin{align}
	&\sum_{j\geq I}m^\tau_{j}+1-\sum_{j\geq I+1}m^\sigma_{j}\leq I,\text{\quad if }\sum_{j\geq I}m^\sigma_{j}<\sum_{j\geq I}m^\tau_{j}+1.\label{cond3}
\end{align}

\medskip

Now, we define the map $f_{n,k,k}$ as follows:

$$f_{n,k,k}(\sigma,\tau)=(\theta,\pi),$$ where :

\medskip

\begin{itemize}
\item $\theta$ is obtained by the following modifications on $\tau$:\\
\begin{itemize}
\item For the entry $(j+1)$, from $(I+2)$ to $n$ in ascending order: if the entry $(j+1)$ from $\sigma$ is overlined, then overlined the entry $(j+1)$ from $\tau$, else not overlined the entry $(j+1)$ from $\tau$.
\item For the entry $(j+1)$, from $n$ to $(I+2)$ in decreasing order: the entry $j+1$ moves $m^\tau_j$ positions to the right and overlined, if the entry $j+1$ is overlined.
Else, the entry $j+1$ only moves $m^\tau_j$ positions to the right.
\item For the entry $j+1=I+1$: the entry $j+1$ moves $\left(\sum_{i=1}^Im^\sigma_i-\sum_{i=1}^Im^\tau_i-1\right)$ positions to the left and overlined, if the entry $j+1$ is overlined. Else, the entry $j+1$ only moves $\left(\sum_{i=1}^Im^\sigma_i-\sum_{i=1}^Im^\tau_i-1\right)$ positions to the left. 
\item For the entry $(j+1)$, from $(I+2)$ to $n$ in ascending order: the entry $j+1$ moves $m_j^\sigma$ positions to the left and overlined, if the entry $j+1$ is overlined.
Else, the entry $j+1$ only moves $m_j^\sigma$ positions to the left.
\end{itemize}

\bigskip

\item $\pi$ is obtained by the following modifications on $\sigma$:
\begin{itemize}
\item For the entry $(j+1)$, from $(I+2)$ to $n$ in ascending order: if the entry $(j+1)$ from $\tau$ is overlined, then overlined the entry $(j+1)$ from $\sigma$, else not overlined the entry $(j+1)$ from $\sigma$.
\item For the entry $(j+1)$, from $n$ to $(I+2)$ in decreasing order: the entry $j+1$ moves $m^\sigma_j$ positions to the right and overlined, if the entry $j+1$ is overlined. Else, the entry $j+1$ only moves $m^\sigma_j$ positions to the right.
\item For the entry $j+1=I+1$: the entry $j+1$ moves $\left(\sum_{i=1}^Im^\sigma_i-\sum_{i=1}^Im^\tau_i-1\right)$ positions to the right and overlined, if the entry $j+1$ is overlined. Else, the entry $j+1$ only moves $\left(\sum_{i=1}^Im^\sigma_i-\sum_{i=1}^Im^\tau_i-1\right)$ positions to the right. 
\item For the entry $(j+1)$, from $(I+2)$ to $n$ in ascending order: the entry $j+1$ moves $m_j^\tau$ positions to the left and overlined, if the entry $j+1$ is overlined.
Else, the entry $j+1$ only moves $m_j^\tau$ positions to the left.
\end{itemize}
\end{itemize}
\bigskip
To prove that $f_{n,k,k}$ is well defined we only need to check that:
\begin{align}
0\leq \sum_{j=1}^Im^\sigma_j-\sum_{j=1}^Im_j^\tau-1\leq m^\sigma_I  \label{eq2}
\end{align}
and
\begin{align}
&m_I^\tau+\sum_{j=1}^Im^\sigma_j-\sum_{j=1}^Im_j^\tau-1\leq I.\label{eq3}
\end{align}
Now, we simplify the inequality \eqref{eq2} as follows:
\begin{align*}
\sum_{j=1}^Im^\sigma_j-\sum_{j=1}^Im_j^\tau-1&=\sum_{j=1}^{n-1}m^\sigma_j-\sum_{j=I+1}^{n-1}m^\sigma_j-\sum_{j=1}^{n-1}m_j^\tau+\sum_{j=I+1}^{n-1}m_j^\tau-1\\&=\left(k+1-\sum_{j=I+1}^{n-1}m^\sigma_j\right)-\left( k-1-\sum_{j=I+1}^{n-1}m_j^\tau\right)-1\\&=1-\sum_{j=I+1}^{n-1}m^\sigma_j+\sum_{j=I+1}^{n-1}m_j^\tau.\label{eq4}
\end{align*}
\medskip
Then inequality \eqref{eq3} becomes:
\begin{align*}
m_I^\tau+\sum_{j=1}^Im^\sigma_j-\sum_{j=1}^Im_j^\tau-1&=m_I^\tau+1-\sum_{j=I+1}^{n-1}m^\sigma_j+\sum_{j=I+1}^{n-1}m_j^\tau\\&=1-\sum_{j=I+1}^{n-1}m^\sigma_j+\sum_{j=I}^{n-1}m_j^\tau. 
\end{align*}
\medskip
Now, we show inequality \eqref{eq2}, i.e.
\begin{align*}
0\leq \sum_{j=1}^Im^\sigma_j-\sum_{j=1}^Im_j^\tau-1\leq m^\sigma_I.
\end{align*}

Here, we have two cases :

\noindent\textbf{First case :} If  $\sum_{j\geq I}m^\sigma_{j}\geq \sum_{j\geq I+1}m^\tau_{j}+1$, we have:
\begin{align*}
1-\sum_{j\geq I}m^\sigma_{j}+\sum_{j\geq I+1}m^\tau_{j}\leq0 &\Leftrightarrow 1-\sum_{j\geq I+1}m^\sigma_{j}-m^\sigma_I+\sum_{j\geq I+1}m^\tau_{j}\leq0\\&\Leftrightarrow 1-\sum_{j\geq I+1}m^\sigma_{j}+\sum_{j\geq I+1}m^\tau_{j}\leq m^\sigma_I\\&\Leftrightarrow \sum_{j=1}^{I}m^\sigma_{j}-\sum_{j=1}^{I}m_{j}^{\tau}-1\leq m_{I}^\sigma.
\end{align*}
\noindent\textbf{Second case :} If  $\sum_{j\geq I}m^\sigma_{j}>\sum_{j\geq I+1}m^\tau_{j}+1$, we have :
\begin{align*}
1-\sum_{j\geq I}m^\sigma_{j}+\sum_{j\geq I+1}m^\tau_{j}<0 &\Leftrightarrow 1-\sum_{j\geq I+1}m^\sigma_{j}-m^\sigma_I+\sum_{j\geq I+1}m^\tau_{j}<0\\&\Leftrightarrow 1-\sum_{j\geq I+1}m^\sigma_{j}+\sum_{j\geq I+1}m^\tau_{j}<m^\sigma_I\\&\Leftrightarrow \sum_{j=1}^{I}m^\sigma_{j}-\sum_{j=1}^{I}m_{j}^{\tau}-1< m_{I}^\sigma.
\end{align*}
It remains to show that :
\begin{align*}
\sum_{j=1}^Im^\sigma_j-\sum_{j=1}^Im_j^\tau-1\geq0.
\end{align*}

For this, we have two cases:

\noindent\textbf{First case :} If  $\sum_{j=I+1}^{n-1}m^\sigma_{j} < \sum_{j=I+2}^{n-1}m^\tau_{j}+1$, we have :
\begin{align*}
-\sum_{j=I+1}^{n-1}m^\sigma_{j} > -\sum_{j=I+2}^{n-1}m^\tau_{j}-1 &\Leftrightarrow 1-\sum_{j=I+1}^{n-1}m^\sigma_{j}>-\sum_{j=I+2}^{n-1}m^\tau_{j}\\&\Leftrightarrow 1-\sum_{j=I+1}^{n-1}m^\sigma_{j}+\sum_{j=I+2}^{n-1}m^\tau_{j}>0 \\&\Leftrightarrow 1-\sum_{j=I+1}^{n-1}m^\sigma_{j}+\sum_{j=I+1}^{n-1}m^\tau_{j}>m_{I+1}^\tau\geq0\\&\Leftrightarrow 1-\sum_{j=I+1}^{n-1}m^\sigma_{j}+\sum_{j=I+1}^{n-1}m^\tau_{j}>0\\&\Leftrightarrow\sum_{j=1}^Im^\sigma_j-\sum_{j=1}^Im_j^\tau-1>0.
\end{align*}
\noindent\textbf{Second case :} If $\sum_{j=I+1}^{n-1}m^\sigma_{j} \geq \sum_{j=I+2}^{n-1}m^\tau_{j}+1$ and $\sum_{j=I+1}^{n-1}m^\sigma_{j} < \sum_{j=I+1}^{n-1}m^\tau_{j}+1$, we have
\begin{align*}
1-\sum_{j=I+1}^{n-1}m^\sigma_{j}+\sum_{j=I+1}^{n-1}m^\tau_{j}>0&\Leftrightarrow\sum_{j=1}^Im^\sigma_j-\sum_{j=1}^Im_j^\tau-1>0.
\end{align*}

\medskip

We show the inequality \eqref{eq3}:\hspace{1cm}$m_I^\tau+\sum_{j=1}^Im^\sigma_j-\sum_{j=1}^Im_j^\tau-1\leq I$.

Here, we have two cases:

\noindent\textbf{First case :} If  $\sum_{j\geq I}m^\sigma_{j}< \sum_{j\geq I}m^\tau_{j}+1$, using the condition \eqref{cond3}, we obtain:
\begin{align*}
\sum_{j\geq I}m^\tau_{j}+1-\sum_{j\geq I+1}m^\sigma_{j}\leq I\Leftrightarrow m_{I}^\tau+\sum_{j=1}^{I}m_{j}^\sigma-\sum_{j=1}^{I}m^\tau_{j}-1\leq I.
\end{align*}
\noindent\textbf{Second case :} If  $\sum_{j\geq I}m^\sigma_{j}\geq\sum_{j\geq I}m^\tau_{j}+1$, we have :
\begin{align*}
1-\sum_{j\geq I}m^\sigma_{j}+\sum_{j\geq I}m^\tau_{j}\leq0 &\Leftrightarrow 1-\sum_{j\geq I+1}m^\sigma_{j}+\sum_{j\geq I}m^\tau_{j}\leq m_{I}^\sigma\leq I\\&\Leftrightarrow 1-\sum_{j\geq I+1}m^\sigma_{j}+\sum_{j\geq I}m^\tau_{j}\leq I \\&\Leftrightarrow m_{I}^\tau+\sum_{j=1}^{I}m^\sigma_{j}-\sum_{j=1}^{I}m^\tau_{j}-1\leq I. 
\end{align*}

\medskip

Now, we must check that $(\theta, \tau)\in I_{B'}(n,k)\times I_{B'}(n,k)$ which means that the number of inversions in $\theta$ is $k$ and the number of inversions in $\pi$ is $k$.

\medskip

The number of inversions in $\theta$ is calculated as follows
\begin{align*}
\sum_{j=1}^{I}m^\tau_j+\left(\sum_{j=1}^Im^\sigma_j-\sum_{j=1}^Im^\tau_j-1\right)+\sum_{j=I+1}^{n-1}m^\sigma_j&=\sum_{j=1}^Im^\sigma_j+\sum_{j=I+1}^{n-1}m^\sigma_j-1\\
&=\sum_{j=1}^{n-1}m^\sigma_j-1\\&=(k+1)-1=k.
\end{align*}

\medskip

The number of inversions of $\pi$ is calculated as follows
\begin{align*}
\sum_{j=1}^{I}m^\sigma_j-\left(\sum_{j=1}^Im^\sigma_j-\sum_{j=1}^Im^\tau_j-1)\right)+\sum_{j=I+1}^{n-1}m^\tau_j&=\sum_{j=1}^{n-1}m^\tau_j+1\\
&=(k-1)+1=k.
\end{align*}

\medskip

Then $(\theta, \tau)\in I_{B'}(n,k)\times I_{B'}(n,k)$.
Furthermore, $f^{-1}_{n,k,k}=f_{n,k+1,k-1}$, which gives us the injectivity of $f_{n,k,k}$ and completes the proof. 
\end{proof}

\medskip

Here, we give an illustrative example.
\begin{example}
	Let $\sigma=\overline{3}\hspace{0.05cm} \overline{2}\hspace{0.05cm}\overline{4}\hspace{0.05cm}\overline{5}1$  a permutation with $5$ overlined inversions and let $\tau =12\overline{5}\hspace{0.05cm} \overline{4}3$  a permutation with $3$ overlined inversions.
	
\medskip
	
	\begin{tabularx}{\linewidth}{@{}l|X@{}}
		$	\begin{aligned}[t]
		&\quad	m_1^\sigma =1\quad \\	&\quad	m_2^\sigma =2\quad \\	&\quad	m_3^\sigma =1\quad\\
		&\quad	m_4^\sigma =1\quad 
		\end{aligned}$
		&
		$	\begin{aligned}[t]
		&\quad	m_1^\tau =0\quad \\	&\quad	m_2^\tau =0\quad \\	&\quad	m_3^\tau =1\quad\\
		&\quad	m_4^\tau =2\quad
		\end{aligned}$
	\end{tabularx}
	
\bigskip

First, we must find the index I, we have :

\bigskip

For $I=4$, $I+1=\overline{5}$, so $m_4^\sigma =1=0+1$ then $I=4$ does not satisfy the second condition \eqref{cond2}.

\medskip

For $I=3$, $I+1=\overline{4}$, so $m_3^\sigma+m_4^\sigma =2<m_4^\tau+1=3$ then $I=3$ does not satisfy the second condition \eqref{cond2}.

\medskip

For $I=2$, $I+1=\overline{2}$, so $m_2^\sigma+m_3^\sigma+m_4^\sigma =4=m_3^\tau+m_4^\tau+1=4$ then $I=2$ does not satisfy the second condition \eqref{cond2}.

\medskip

For $I=1$, $I+1=\overline{3}$, so $m_1^\sigma+m_2^\sigma+m_3^\sigma+m_4^\sigma =5>m_2^\tau+m_3^\tau+m_4^\tau+1=4$ then $I=1$ satisfies the second condition \eqref{cond2}. Thus, $I = 1$ and our application gives :

\bigskip

\begin{center}
$\overline{3}\hspace{0.05cm} \overline{2}\hspace{0.05cm}\overline{4}\hspace{0.05cm}\overline{5}1 \to$ $3\overline{2}\hspace{0.05cm}\overline{4}\hspace{0.05cm}\overline{5}1\to$ 
$3\overline{2}\hspace{0.05cm}\overline{4}1\hspace{0.05cm}\overline{5}\to$
$3\overline{2}1\overline{4}\hspace{0.05cm}\overline{5}\to$
$\overline{2}13\overline{4}\hspace{0.05cm}\overline{5}\to$
$\overline{2}1\overline{4}3\overline{5}\to$
$\overline{2}1\overline{5}\hspace{0.05cm}\overline{4}3$
\end{center}

\bigskip

\begin{center}
$12\overline{5}\hspace{0.05cm} \overline{4}3\to$
 $12\overline{5}\hspace{0.05cm}\overline{4}\hspace{0.05cm}\overline{3}\to$ 
$12\overline{4}\hspace{0.05cm}\overline{3}\hspace{0.05cm}\overline{5}\to$
$12\overline{3}\hspace{0.05cm}\overline{4}\hspace{0.05cm}\overline{5}\to$
$\overline{3}12\overline{4}\hspace{0.05cm}\overline{5}\to$
$\overline{3}1\overline{4}2\overline{5}\to$
$\overline{3}1\overline{4}\hspace{0.05cm}\overline{5}2$
\end{center}

\bigskip

Hence, $f_{5,4,4}(\overline{3}\hspace{0.05cm} \overline{2}\hspace{0.05cm}\overline{4}\hspace{0.05cm}\overline{5}1,12\overline{5}\hspace{0.05cm} \overline{4}3)=(\overline{2}1\overline{5}\hspace{0.05cm}\overline{4}3,\overline{3}1\overline{4}\hspace{0.05cm}\overline{5}2)$,
where $\theta=\overline{2}1\overline{5}\hspace{0.05cm}\overline{4}3$ has $4$ overlined inversions, and $\pi=\overline{3}1\overline{4}\hspace{0.05cm}\overline{5}2$ has $4$ overlined inversions.	
\end{example}

\medskip

Since a log-concave sequence is unimodal \cite{bren2}, the previous theorem immediately leads to the following corollary:
\begin{corollary}
The sequence of over-Mahonian numbers $\{i_{B'}(n,k)\}_{k}$ is unimodal in $k$.
\end{corollary}

\begin{remark}
Since there is a bijection between the permutations with overlined inversions and the paths as shown in Subsection \ref{path}, then the previous theorem can be proved using Sagan's paths approach \cite{sag}  by adding a new condition to his involution as follows. (See also the approaches proposed in \cite{GA1,BZB}).

\bigskip

Let $P_1,P_2\in P_{n,k}^{B'}$. And let $u \overset{{P}}{\longrightarrow} v$ denote  $P$ has initial vertex $u$ and final vertex $v$.
\begin{definition}\label{df1}
	Given $u_1 \overset{{P_1}}{\longrightarrow} v_1$ and $u_2 \overset{{P_2}}{\longrightarrow} v_2$. Then define the involution $\varphi_I(P_1,P_2)=(P^{'}_{1},P^{'}_{2})$
	where :
	\begin{enumerate}
\item[] $$P^{'}_{1}=u_1 \overset{{P_1}}{\longrightarrow} v_0 \overset{{P_2}}{\longrightarrow} v_2 \text{\ and\ }P^{'}_{2}= u_2 \overset{{P_2}}{\longrightarrow} v_0 \overset{{P_1}}{\longrightarrow} v_1,$$
i.e., switches the portions of
		$P_1$ and $P_2$ after $v_0$ (see, Figure \ref{Fig2}), where $v_0$ is the last vertex of $P_1 \cap P_2$.	
	\end{enumerate}
	
\end{definition}
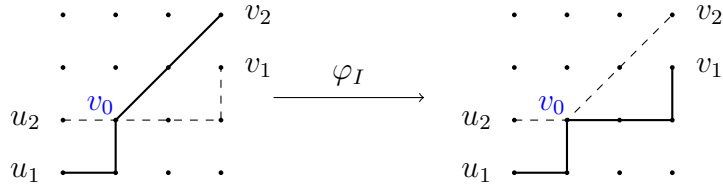
\begin{figure}[h!]
		\begin{center}
			
			\begin{tikzpicture}
			
			\draw[dashed]  (-0.7,0.7)--(0,0.7)--(0.7,0.7)--(1.4,0.7)--(1.4,1.4);
			\draw [line width=0.8pt] (-0.7,0)--(0,0)--(0,0.7)--(0.7,1.4)--(1.4,2.1);
			
			\fill[] (-0.7,0) circle (0.9pt);
			\fill[] (-0.7,0.7) circle (0.9pt);
			\fill[] (-0.7,1.4) circle (0.9pt);
			\fill[] (-0.7,2.1) circle (0.9pt);

			\fill[] (0,0) circle (0.9pt);
			\fill[] (0,0.7) circle (0.9pt);
			\fill[] (0,1.4) circle (0.9pt);
			\fill[] (0,2.1) circle (0.9pt);

			\fill[] (0.7,0) circle (0.9pt);
			\fill[] (0.7,0.7) circle (0.9pt);
			\fill[] (0.7,1.4) circle (0.9pt);
			\fill[] (0.7,2.1) circle (0.9pt);
			\fill[] (1.4,0) circle (0.9pt);
			\fill[] (1.4,0.7) circle (0.9pt);
			\fill[] (1.4,1.4) circle (0.9pt);
			\fill[] (1.4,2.1) circle (0.9pt);


			\node at (-1.2,0){$u_1$};
			\node at (-1.2,0.7){$u_2$};
			\node at (1.9,1.4){$v_1$};
			\node at (1.9,2.1){$v_2$};
			
			\node at (-0.2,0.9){$\textcolor[rgb]{0.00,0.00,1.00}{v_0}$};

			\begin{scope}[xshift=1.9cm]
			\draw [->] (0.2,1)--(2.2,1);
			\node[rectangle] at (1.2,1.3) {$\varphi_I$};
			
			\begin{scope}[xshift=4.1cm]
			\draw[dashed]    (-0.7,0.7)--(0,0.7)--(0.7,1.4)--(0.7,1.4)--(1.4,2.1);
			\draw [line width=0.8pt] (-0.7,0)--(0,0)--(0,0.7)--(0.7,0.7)--(1.4,0.7)--(1.4,1.4);
			
			\fill[] (-0.7,0) circle (0.9pt);
			\fill[] (-0.7,0.7) circle (0.9pt);
			\fill[] (-0.7,1.4) circle (0.9pt);
			\fill[] (-0.7,2.1) circle (0.9pt);
			\fill[] (0,0) circle (0.9pt);
			\fill[] (0,0.7) circle (0.9pt);
			\fill[] (0,1.4) circle (0.9pt);
			\fill[] (0,2.1) circle (0.9pt);
			\fill[] (0.7,0) circle (0.9pt);
			\fill[] (0.7,0.7) circle (0.9pt);
			\fill[] (0.7,1.4) circle (0.9pt);
			\fill[] (0.7,2.1) circle (0.9pt);
			
			\fill[] (1.4,0) circle (0.9pt);
			\fill[] (1.4,0.7) circle (0.9pt);
			\fill[] (1.4,1.4) circle (0.9pt);
			\fill[] (1.4,2.1) circle (0.9pt);

			\node at (-1.2,-0){$u_1$};
			\node at (-1.2,0.7){$u_2$};
		    \node at (1.9,1.4){$v_1$};
			\node at (1.9,2.1){$v_2$};

			\node at (-0.2,0.9){$\textcolor[rgb]{0.00,0.00,1.00}{v_0}$};

			\end{scope}
			\end{scope}
			\end{tikzpicture}	
		\end{center}
		
		\caption{The involution $\varphi_{I}$.}\label{Fig2}
	\end{figure}
\end{remark} 
\section{The remark and the question about the mode}
In the previous section, we have established the unimodality and log-concavity properties of the sequence of over-Mahonian numbers. However the number and location of the
modes of this sequence remains a question to be answered. Generally, it is not easy to find the number and location of modes.

\medskip

This lets us to finish this paper by the following question.

\smallskip

\noindent{\bf Question.} Find the number and location of modes of the unimodal sequence $\left\{ i_{B'}(n,k)\right\}_k$.
\section*{Acknowledgment}
The authors would like to thank the referees for many valuable remarks and suggestions to improve the original manuscript. This work was supported by DG-RSDT (Algeria), PRFU Project C00L03UN180120220002 and PRFU Project C00L03UN190120230012.

\end{document}